\documentclass[12pt]{article}

\usepackage{amsmath,amsthm,amssymb,graphicx}

\def\N{\mathbb N}
\def\N0{\mathbb N_0}

\pdfoutput=1

\title{\bf An Extended Pruess Method   \\
for Sturm-Liouville Problems}
\author{Robert Carlson \\
Department of Mathematics \\ 
University of Colorado at Colorado Springs \\
carlson@math.uccs.edu}

\theoremstyle{definition}

\theoremstyle{remark}

 \numberwithin{equation}{section}

\begin{document}

\maketitle

\begin{abstract}

A new version of the piecewise approximation (Pruess) method is developed for calculating
eigenvalues of Sturm-Liouville problems.  The usual piecewise constant or piecewise linear
potential approximations are replaced by translates of $2/\cos ^2(x)$, whose corresponding
eigenvalue equation has elementary solutions.

\end{abstract}

\vskip 50pt


\vskip 10pt


\newpage

\section{Introduction}

The basic Sturm-Liouville problem asks for the eigenvalues and eigenfunctions of the second order linear equation
\begin{equation} \label{eqA}
-y'' + p(x)y = \lambda y, \quad 0 \le x \le 1,
\end{equation}
satisfying the separated boundary conditions  
\begin{equation} \label{bcA}
a_0y(0) + a_1y'(0) = 0, \quad b_0y(1) + b_1y'(1) = 0.
\end{equation}
The constants $a_j, b_j$ should be real, and not both zero at either endpoint.
The potential $p(x)$ should be real valued; although singular potentials are interesting, our examples will be continuous.
A closely related problem asks for periodic eigenfunctions and the associated eigenvalues, with the 
boundary conditions
\begin{equation} \label{bcB}
y(0) = y(1), \quad  y'(0) = y'(1).
\end{equation}
The book \cite{Pryce} provides a systematic account of numerical methods for solving Sturm-Liouville problems. 

Popular algorithms known as piecewise constant midpoint methods, or Pruess methods, are based on a piecewise 
constant approximation $\tilde{p}(x)$ of the coefficient $p(x)$.
This algorithmic theme was introduced in \cite{Canosa} and eventually developed into well analyzed packages \cite{Ledoux09, Ledoux10, Marletta92, Pruess93}.
When $\tilde p$ is piecewise constant, the solutions of \eqref{eqA} have a simple piecewise description in terms of the trigonometric functions
$\cos (x), \sin(x)$.   Taking advantage of the simple explicit formulas for the basis,
a shooting method can provide an efficient algorithm for accurately computing eigenvalues and eigenfunctions.

Piecewise constant approximations are rather crude, so it is not surprising that more refined approximations were considered.
Piecewise polynomial approximations were analyzed in \cite{Pruess73}, but the convenient basis of solutions is generally absent, making corresponding
algorithms much less attractive.  These difficulties are reduced for piecewise linear approximation \cite{Ledoux67,Ledoux06} 
because the bases can be expressed efficiently in terms of Airy functions. 
 
This work considers a different type of piecewise approximation, focusing on the existence of an elementary basis of solutions to \eqref{eqA}.
By using commutation methods, also known as Darboux transformations, it is possible to find examples of \eqref{eqA} with nonconstant
coefficients whose solutions are elementary functions for all values of $\lambda $.   
A well-known example is 
\begin{equation} \label{Bessform}
-y'' + \frac{2}{x^2}y = \lambda y, 
\end{equation}
which has solutions
\begin{equation} \label{Bsoln}
y_1(x,\lambda ) = \sin (\sqrt{\lambda } x) + \frac{\cos (\sqrt{\lambda } x)}{\sqrt{\lambda } x} , \quad 
y_2(x,\lambda ) = \cos (\sqrt{\lambda } x) - \frac{\sin (\sqrt{\lambda} x)}{\sqrt{\lambda } x} .
\end{equation}
The existence of such examples suggests that certain nonpolynomial piecewise approximations of $p(x)$ in \eqref{eqA}
might provide an extension of the piecewise constant Pruess method with better approximation of $p(x)$,
elementary expressions for the solutions, and more efficient computation of the eigenvalues.

This work describes such an algorithm and reports on its performance in a number of numerical experiments.
The piecewise approximations of the extended algorithm have the form 
\begin{equation} \label{localmodel}
p_k(x) = \alpha _k + \frac{2}{\cos ^2(x - \xi _k)}, \quad x_k \le x < x_{k+1},
\end{equation} 
taking advantage of the fact that the equation
\begin{equation} \label{basecase}
- y'' + \frac{2}{\cos ^2(x)}y = \lambda y 
\end{equation}
has a basis of solutions $Y_1,Y_2$ with
\begin{equation} \label{thebasis}
Y_1(x,\lambda )  =  \cos(\sqrt{\lambda }x) + \tan(x)\frac{\sin(\sqrt{\lambda }x)}{\sqrt{\lambda }},
\end{equation}
\[Y_2(x,\lambda ) = \frac{1}{1 - \lambda} \tan (x) \cos (\sqrt{\lambda }x) - \frac{\lambda }{1 - \lambda } \frac{\sin(\sqrt{\lambda }x)}{\sqrt{\lambda }} ,\]
\[Y_1'(x,\lambda )  =  \tan(x)\cos(\sqrt{\lambda }x) + [\frac{1}{\cos^2(x)} - \lambda ] \frac{\sin(\sqrt{\lambda }x)}{\sqrt{\lambda }}],\] 
\[Y_2'(x,\lambda ) = [\frac{1}{\cos ^2(x)} - \lambda ]\frac{\cos (\sqrt{\lambda }x)}{1-\lambda } - \frac{\lambda }{1 - \lambda} \tan (x) \frac{\sin (\sqrt{\lambda }x)}{\sqrt{\lambda }} .\]
On the interval $(-\pi /2,\pi/2)$ the derivative values of $2/\cos ^2(x)$ increase monotonically from $-\infty $ to $\infty $.
For each subinterval $[x_k,x_{k+1}]$ our algorithm computes a derivative estimate for $p(x)$, then selects $\xi _k$ in \eqref{localmodel}
to match the estimate.  The constant $\alpha _k$ is then adjusted to optimize the nearly linear local approximation.

The subsequent sections describe the background ideas, give a description of the algorithm, and report on its performance.
After a review of piecewise approximations of \eqref{eqA} and the corresponding transition matrices,
the general commutation method is considered.   Examples are provided.  When the examples are more
complex than \eqref{basecase} the formulas for solutions to the eigenvalue equation \eqref{eqA} quickly grow in complexity;
the properties of the corresponding potentials are hard to discern.  The extended algorithm presented here uses \eqref{basecase};
it produces a piecewise approximation of $p(x)$ which is nearly linear, but with solutions which are elementary functions instead of Airy functions.

Algorithm accuracy was compared for several variants of both the piecewise constant Pruess method and the new algorithm.
For both methods the sample points were chosen in two ways: uniformly spaced, and adaptively chosen to optimize piecewise approximation.   
The lowest $25$ eigenvalues were computed for five regular Sturm-Liouville problems.
The reported data comprises eigenvalues $1,2,3,12$ and $25$.  The lowest eigenvalues typically show the largest errors.
While the new algorithm does quite well in approximation of the potentials, all of the algorithms had similar performance for computation
of the eigenvalues.  As long as eigenvalue computation is the goal, the simplicity and speed of the piecewise constant Pruess method
make it preferable.
  
\section{The basic idea}

Suppose the interval $[0,1]$ is partitioned into $K$ subintervals with endpoints $0 = x_0 < x_1 < \dots < x_K = 1$.
Assume that the function $p(x)$ from \eqref{eqA} has a piecewise definition based on this partition,
\[p(x) = p_k(x), \quad x_k \le x < x_{k+1}, \quad k = 0,\dots , K-1.\]
The functions $p_k(x)$ are assumed to have a continuous extension to the closed interval $[x_k,x_{k+1}]$.
On each subinterval there is a basis $y_k(x,\lambda ), z_k(x,\lambda )$ of solutions to \eqref{eqA} satisfying the initial conditions 
\[y_k(x_k, \lambda ) = 1, \quad y_k'(x_k,\lambda ) = 0,\]
\[z_k(x_k, \lambda ) = 0, \quad z_k'(x_k, \lambda ) = 0.\]
It will be convenient to package these functions into a matrix
\begin{equation}  \label{Mono}
\begin{pmatrix} y_k(x,\lambda ) & z_k(x,\lambda ) \cr y_k'(x,\lambda ) & z_k'(x,\lambda )  \end{pmatrix} 
\end{equation}
Because $p_0(x)$ has a continuous extension to the closed interval $[x_0,x_1]$,
the solutions $y_0,z_0$ extend to $x_1$, prescribing initial data 
\[\begin{pmatrix} y(x_1,\lambda ) & z(x_1,\lambda ) \cr y'(x_1,\lambda ) & z'(x_1,\lambda )  \end{pmatrix} 
= \begin{pmatrix} y_0(x_1,\lambda ) & z_0(x_1,\lambda ) \cr y_0'(x_1,\lambda ) & z_0' (x_1,\lambda ) \end{pmatrix} ,\]
for the continuation to the next interval, where the solution matching values and derivatives at $x_1$ takes the form
\[ \begin{pmatrix} y(x,\lambda ) & z(x,\lambda ) \end{pmatrix} = \begin{pmatrix} y_1(x,\lambda ) & z_1(x,\lambda ) \end{pmatrix} 
\begin{pmatrix} y_0(x_1,\lambda ) & z_0(x_1,\lambda ) \cr y_0'(x_1,\lambda ) & z_0' (x_1,\lambda ) \end{pmatrix} .\]
By induction we find that the values and derivatives of $y$ and $z$ at $1$ are given by the product of matrices, 
 \[ \begin{pmatrix} y(1,\lambda ) & z(1,\lambda ) \cr y'(1,\lambda ) & z'(1,\lambda ) \end{pmatrix} = \prod _{k=0}^{K-1} T_k(\lambda ),
 \quad T_k(\lambda ) = \begin{pmatrix} y_{k}(x_{k+1},\lambda ) & z_{k}(x_{k+1},\lambda ) \cr y_{k}'(x_{k+1},\lambda ) & z_{k}' (x_{k+1},\lambda ) \end{pmatrix},\]
where the factors in the product have indices decreasing from left to right, 
\begin{equation} \label{Tdef}
 \prod _{k=0}^{K-1} T_k(\lambda ) = T_{K-1}(\lambda ) \cdots T_1(\lambda )T_0(\lambda ), 
 \end{equation}

The piecewise constant approximation methods use $p_k(x) = \alpha _k $.
If $l_k = x_{k+1}-x_k$, the matrices $T_k$ then have the elementary form
\begin{equation} \label{pcT}
T_k(\lambda ) = \begin{pmatrix} \cos (l_k \sqrt{\lambda - \alpha _k}  ) & \sin (l_k \sqrt{\lambda - \alpha _k} ) /\sqrt{\lambda - \alpha _k}
 \cr - \sqrt{\lambda - \alpha _k}  \sin (l_k \sqrt{\lambda - \alpha _k}) & \cos (l_k \sqrt{\lambda - \alpha _k}) \end{pmatrix}.
 \end{equation}
The next section will discuss the construction of model potentials leading to roughly linear piecewise approximations of the potential $p(x)$
which retain elementary (but more complex) formulas for the matrices $T_k(\lambda)$.

\section{Commutation methods}

It is possible to generate a large collection of differential operators whose  eigenvalue equations  \eqref{eqA} have
elementary solutions like \eqref{thebasis}.  These operators and their solutions are generated from the case $p(x) = 0$ using 
the commutation method, also known as a Darboux transformation.
Such ideas have been used repeatedly in the past \cite{BC}, \cite{Carlson82}, \cite[p. 425-444]{Gesz1},  \cite[p. 88-91]{PT}.

\subsection{General idea}

Suppose $L_0$ is a second order operator of the form $L_0 = -D^2 + p_0$, and $y_1$ is a nontrivial solution of
\[(L_0 - \mu _1 ) y = 0.\]
Then $L_0-\mu _1$ may be factored,
\[L_0-\mu _1  = F_-F_+  , \quad F_{\pm} = \pm D - \frac{y_1'}{y_1}, \quad p_0 - \mu _1 = (\frac{y_1'}{y_1})' + (\frac{y_1'}{y_1})^2.\]

If $z$ satisfies a second eigenvalue equation $L_0z = \lambda z$, then $(F_-F_+ + \mu _1 )z = \lambda z$.
Applying $F_+$ to both sides gives 
\[[F_+F_- + \mu _1] F_+ z = \lambda F_+z ,\]
so that the new operator $L_1 = -D^2 + p_1$ has a solution  $F_+z$ for the eigenvalue equation $L_1F_+z = \lambda F_+z$, where
\begin{equation} \label{newop}
L_1  =  F_+F_- + \mu _1 =  -D^2 - (\frac{y_1'}{y_1})' + (\frac{y_1'}{y_1})^2 + \mu _1 = -D^2 - p_0 + 2\mu _1 + 2(\frac{y_1'}{y_1})^2.
\end{equation}
That is 
\[p_1 = -p_0 + 2\mu _1 + 2(\frac{y_1'}{y_1})^2.\]
Notice too that if $\lambda \not= \mu _1$, then the linear map $z \to F_+z$ is invertible since 
\[\begin{pmatrix}F_+z \cr (F_+z)' \end{pmatrix} = A \begin{pmatrix} z \cr z' \end{pmatrix},
\quad A = \begin{pmatrix} -y_1'/y_1 & 1 \cr p - \lambda -(y_1'/y_1)' & -y_1'/y_1 \end{pmatrix},\]
with $\det (A) = \lambda - \mu _1 $. 

If $\mu _2 \not= \mu _1 $ the process may be repeated.  Adding parameters for clarity, the mapping on $z$ is given by
\[F_+(\mu _1, \mu _2)z(\lambda ) = \Bigl [ D - \frac{(F_+(\mu _1)y_2(\mu _2 ))'}{F_+(\mu _1) y_2(\mu _2)} \Bigr ] F_+(\mu _1 ) z .\] 
The operator $F_+(\mu _1,\mu _2)$, which is second order with leading coefficient $1$, annihilates $y_1(\mu _1)$ and $y_2(\mu _2)$, so
\begin{equation} \label{doublecom}
F_+(\mu _1, \mu _2) z = \frac{1}{y_1(\mu _1 )y_2'(\mu _2 ) -  y_1'(\mu _1)y_2(\mu _2 )}\det \begin{pmatrix} z & y_1 & y_2 \cr 
z' & y_1' & y_2' \cr
z'' &  y_1'' & y_2''  
\end{pmatrix} .
\end{equation}
Notice that $F_+(\mu _1,\mu _2) = F_+(\mu _2 ,\mu _1)$, and this construction extends to maps $F_+(\mu _1,\dots ,\mu _n)$,
with
\begin{equation} \label{ncom}
F_+(\mu _1, \dots ,\mu _n) z 
\end{equation}
\[= (-1)^n
\det \begin{pmatrix}
y_1 & \dots & y_n \cr
 \vdots & \dots & \vdots  \cr
y_1^{(n-1)} & \dots & y_n^{(n-1)}  
\end{pmatrix} ^{-1}
\det \begin{pmatrix} z & y_1 & \dots & y_n \cr
z^{(1)} & y_1^{(1)} & \dots & y_n^{(1)} \cr 
\vdots & \vdots & \dots & \vdots \cr
z^{(n)} & y_1^{(n)} & \dots & y_n^{(n)}  
\end{pmatrix} .\]

The operators $-D^2 + p_n$ resulting from an $n$-th commutation of the initial operator $-D^2+p_0$ have a recursive description.
For the $n$-th commutation, a solution of $(-D^2 + p_{n-1})y = \mu _ny$ has the form
$z_n = F_+(\mu _1,\dots , \mu_{n-1})y_n$ where $(-D^2 + p_0)y_n= \mu _ny_n$.  
From \eqref{newop} we see that the 
$n$-th commutation results in a new potential
\begin{equation} \label{potn}
p_n = - (\frac{z_n'}{z_n})' + (\frac{z_n'}{z_n})^2 + \mu _n
=  -p_{n-1} + 2\mu _n + 2(\frac{z_n'}{z_n})^2.
\end{equation}

Suppose the coefficient $p_0$ is continuous on $[\alpha , \beta ]$.  
If $y_1(\mu _1)$ has a zero on this interval then $F_+(\mu _1)$ and the modified operator \eqref{newop} will have singularities.
The situation is different for the double commutation \eqref{doublecom}.
Notice that $F_+(\mu _1,\mu _2)$ is defined when $\mu _1 = \mu _2$ if $y_1(\mu _1)$ and $y_2(\mu _2)$ are linearly independent.
In this case $F_+(\mu _1,\mu _1 )$ agrees with $L - \mu _1$ and the operator is not altered.  
The transformed operator will remain nonsingular on $[\alpha ,\beta ]$ if $y_2(\mu _2)$ is sufficiently close to $y_2(\mu _1)$.

\subsection{Elementary single commutation}
If $p_0 = 0$ and $y_1 = \cos (\sqrt{\mu _1 } x)$, then 
\[F_+ \frac{\sin(\sqrt{\lambda }x)}{\sqrt{\lambda } } =  \cos (\sqrt{\lambda } x) + \sqrt{\mu _1} \tan (\sqrt{\mu _1}x) \frac{\sin(\sqrt{\lambda }x)}{ \sqrt{\lambda }}\]
and 
\[F_+ \cos(\sqrt{\lambda }x) = - \lambda \frac{ \sin (\sqrt{\lambda } x)}{\sqrt{\lambda }} + \sqrt{\mu _1} \tan (\sqrt{\mu _1}x)\cos(\sqrt{\lambda }x)\]
are solutions of 
\begin{equation} \label{com1}
 [-D^2 + \frac{2\mu _1}{\cos ^2(\sqrt{\mu _1}x)}] Y = \lambda Y , \quad p_1 = \frac{2\mu _1}{\cos ^2(\sqrt{\mu _1}x)}.
 \end{equation}
 
 Assuming that $\sqrt{\mu _1} > 0$, the derivative of $p_1(x)$ increases monotonically from $-\infty $ to $\infty $  on the interval 
 $(-\frac{\pi}{ 2 \sqrt{\mu _1}}, \frac{\pi}{2\sqrt{\mu _1}})$.
Simple calculations give
\[ \frac{d}{dx} \frac{2\mu _1}{\cos ^2(\sqrt{\mu _1}x)}  = 4\mu _1^{3/2} \tan(\sqrt{\mu _1} x) [1 + \tan ^2(\sqrt{\mu _1} x) ]\]
and
\[ \frac{d^2}{dx^2} \frac{2\mu _1}{\cos ^2(\sqrt{\mu _1}x)}  
= 4\mu _1^{2}\sec ^2(\sqrt{\mu _1}x)  [1 + 3\tan ^2(\sqrt{\mu _1} x) ] \]
\[ = 4\mu _1^{2}[1 + \tan ^2(\sqrt{\mu _1}x)]  [1 + 3\tan ^2(\sqrt{\mu _1} x) ] \]
Notice that the second derivative is always positive, so a richer family of examples would be needed to locally match 
first and second derivatives of a potential $p(x)$.

In case $\sqrt{\mu _1} = 1$, a basis for $L_1y = \lambda y$ is given by \eqref{thebasis}.
If $y_1 = \sin(\sqrt{\mu _1}x)$ instead of $y_1 = \cos(\sqrt{\mu _1}x)$, then
\[\frac{y_1'}{y_1} = \frac{\sqrt{\mu _1}\cos(\sqrt{\mu_1}x) }{\sin(\sqrt{\mu _1}x)}, \quad 
p_1 = 2\mu _1 + 2(\frac{y_1'}{y_1})^2 = \frac{2\mu _1}{\sin ^2(\sqrt{\mu _1} x)} .\]
Notice that $p_1 = 2/x^2$ when $\mu _1 = 0$.  Repeated commutations can generate potentials
$q_m(x) = \frac{m(m+1)}{x^2}$ for $m = 1,2,3,\dots $.  These cases are closely related to Bessel's equation
of half integer order $\nu = (2m+1)/2$; see  \cite[p.364-365]{WW} for some explicit formulas.

\subsection{Elementary double commutation}

Since single commutation leads to potentials with constrained second derivatives, it is natural to consider examples generated by repeated commutations.
The formulas quickly become unwieldy.  A double commutation can be based on 
\[z_2 =  \sqrt{\mu _2} \cos (\sqrt{\mu _2 } x) + \sqrt{\mu _1} \tan (\sqrt{\mu _1}x) \sin(\sqrt{\mu _2 }x)\]
with
\[z_2' = -\mu _2 \sin (\sqrt{\mu _2}x) \]
\[+ \mu _1 \sec ^2(\sqrt{\mu _1 x}) \sin (\sqrt{\mu _2} x)
+ \sqrt{\mu _1}\sqrt{\mu _2} \tan (\sqrt{\mu _1}x)\cos (\sqrt{\mu _2}x).\]
Transformation of \eqref{com1} results in the new operators $L_2 = -D^2 + p_2$ with
\begin{equation} \label{p2}
p_2 = \frac{-2\mu _1}{\cos ^2(\sqrt{\mu _1}x)} + 2 \mu _2 
\end{equation}
\[ + 2  \tan ^2(\sqrt{\mu _2}x)\Bigl [ \frac{(\mu _1 -\mu _2)  + \mu _1 \tan ^2(\sqrt{\mu _1}x)  
+ \sqrt{\mu _1}\sqrt{\mu _2} \tan (\sqrt{\mu _1}x)\cot (\sqrt{\mu _2}x)
}{\sqrt{\mu _2}  + \sqrt{\mu _1} \tan (\sqrt{\mu _1}x) \tan(\sqrt{\mu _2 }x)} \Bigr ] ^2 .\]

\section{Algorithm description}

\subsection{Table development}

The computations are based on the potential 
\[q(x)  = \frac{2}{\cos ^2(x)},\quad q'(x) = \frac{4\sin(x)}{\cos^3 (x)} = 4\tan (x)[1 + \tan ^2(x)] .\]
On a subinterval interval $[x_k,x_{k+1}] \subset [0,1]$ with midpoint $m_k = (x_k+x_{k+1})/2$
the full model potential will have the form
\[\alpha _k + \frac{2}{\cos^2(x - m_k + z_k)}. \]
The point $z_k$ is chosen so that $q'(z_k)$ approximates $p'(m_k)$.
The selection of $z_k$ is  based on a table of 201 samples $((q')^{1/3}(t_j), t_j)$,
with $q'(t_j)$ in the range $[-1000,1000]$.  The cube root was selected to provide adequate sampling
for moderate values of $q'$.  

A fixed number $N$ of sample points $0 \le x_k \le 1$ were selected.  Some trials used uniformly spaced samples for both
the piecewise constant Pruess method and the extended method. Another set of samples was selected to minimize
\begin{equation} \label{penaltyp}
\sum_{k=0}^{N-1} \int_{x_k}^{x_{k+1}} (p(x) - p(m_k))^2 \ dx 
\end{equation}
for the Pruess method.  For the new method a  
set of samples was selected based on a piecewise linear approximation, the samples chosen to minimize
\begin{equation} \label{penaltyX}
\sum_{k=0}^{N-1} \int_{x_k}^{x_{k+1}} (p(x) - [P_k'(x-m_k) + p(m_k)])^2 \ dx , \quad 
P'_k = \frac{p(x_{k+1}) - p(x_k)}{x_{k+1} - x_k}.
\end{equation}

After the selection of sample points $x_k$, the values $P_k'$ are also used to select the points $z_k$ where $q'(z_k) \simeq p'(m_k)$.
The points $z_k$ are computed by combining linear interpolation with
the table values discussed above. The constants  
\[\alpha _k = p(m_k) - 2[\tan (m_k + z_k) - \tan (z_k - m_k)]/L_k .\] 
are chosen so that the model potential on $[x_k,x_{k+1}]$ 
has the same integral as given by the midpoint rule.

Having selected $x_k,z_k,\alpha _k$, the equation which locally approximates 
\[- y'' + p(x)y = \lambda y\]
on the interval $[x_k,x_{k+1}]$ is
\[-y'' + [\frac{2}{\cos ^2(x - m_k + z_k)} + \alpha _k ]y = \lambda y,\]
or
\[-y'' + \frac{2}{\cos ^2(x- m_k + z_k)} y = \sigma y, \quad \sigma = \lambda - \alpha _k.\] 
Using \eqref{thebasis}, this equation has a basis of solutions $Y_1(x - m_k +z_k,\sigma ), Y_2(x - m_k +z_k,\sigma )$.
After converting to the basis $(W_1,W_2)$ satisfying
\[\begin{pmatrix} W_1 & W_2 \cr W_1' & W_2' \end{pmatrix}(x_k,\lambda ) = \begin{pmatrix} 1 & 0 \cr 0 & 1 \end{pmatrix}.\]
and letting $ l_k = x_{k+1}-x_k$, the transition matrix $T_{k}$ will be
\[ T_k(\lambda ) = \begin{pmatrix} Y_1 & Y_2 \cr Y_1' & Y_2 \end{pmatrix} (x_{k+1} - m_k + z_k, \sigma )  
\begin{pmatrix} Y_2' & -Y_2 \cr -Y_1' & Y_1 \end{pmatrix}(z_k - \frac{l_k}{2},\sigma ).\]

Although the evaluation of $T_k(\lambda )$ in the extended method is not computationally expensive, the evaluation \eqref{pcT} for piecewise constant
approximations is more efficient.  The total algorithm ran about three times as fast for the $16$ subinterval cases if the piecewise constant
approximation was used instead of the extended algorithm.

Trials were run for regular Sturm-Liouville problems on the interval $[0,1]$ using the Dirichlet boundary conditions $y(0) = 0 = y(1)$.
Extensive theoretical studies \cite{PT} of these problems are available.
The eigenvalues are simple, and typically well separated.  
The algorithm used initial samples of the eigenvalue parameter expected to be considerable closer than the actual eigenvalues.
A bisection method searched for roots of $W_2(1,\lambda )$, giving the Dirichlet eigenvalues. 

\section{Performance}

The following five potentials were tested: 
\[(1) \  p(x) = \frac{\pi ^2}{(\pi x+0.1)^2},\]
\[(2) \ p(x) =  1 + \cos (\pi t) + 5\cos(2\pi t) -2 \cos (3\pi t) - 3\cos (4\pi t), \quad t = x - 0.5,\]
\[(3) \ p(x) = \Bigl \{ \begin{matrix} x \sin (1/x), & x \ge 10^{-6}, \cr
0, & x < 10^{-6} \end{matrix} \Bigr\} ,\]
\[(4) \ p(x) = \frac{1}{\cos(x)^2}, \quad (5) \ p(x) = \frac{1}{.2 + \sqrt{x(1-x)}} .\]
In all cases, $0 \le x \le 1$ and the Dirichlet boundary conditions $y(0) = 0 = y(1)$ were used.

The first $25$ eigenvalues were computed.  The main experiments compared the performance of several algorithms 
when the unit interval was partitioned with $16$ subintervals.  Performance was assessed by comparing the computed eigenvalues with 
the results obtained from the Pruess method (U-P128) using $128$ subintervals of equal length.  
Case 1 was used to test the accuracy of U-P128.  This case, considered to be challenging, 
is a rescaling to $[0,1]$ of a problem in Appendix A of \cite{Pryce}.
The rescaling means the numbers calculated here should agree with \cite{Pryce} after division by $\pi ^2$.
A comparison showed a maximum difference of $4$ in the least significant (fifth) digit of the second eigenvalue.
The other eigenvalues shown here had a difference of $1$ in the least significant digit. 

The tables and figures below show a sampling of results.  The figures show plots of the potential $p(x)$ and either the
$16$ equal length subinterval piecewise constant approximation, or the  
extended method with $16$ subintervals, with sample points chosen to minimize the penalty of \eqref{penaltyX}.
In each case there are six trials.  The trials labelled U-P 16, U-P 32, and U-P 128
used the piecewise constant Pruess method with $16$, $32$, and $128$ equal length subintervals respectively.
The trial labelled A-P 16 used 16 subintervals, with sample points chosen to minimize the penalty of \eqref{penaltyp}.
The trial U-X 16 used the extended method with $16$ equal length subintervals.
The trial A-X 16 used the extended method with $16$ subintervals, with sample points chosen to minimize the penalty of \eqref{penaltyX}.

The tables show the values of $5$ eigenvalues: the first three, number $12$ and number $25$.  
For these methods the smaller eigenvalues are more challenging than the larger ones.
The U-P 32 algorithm is included as an indication of the improvement when the number of equal length subintervals is doubled.
The bottom of the table indicates which of the $16$ subinterval algorithms was most accurate.
Overall the performances were similar.

The extended algorithm generally did a good job providing a piecewise approximation of the potential.
If such an approximation were the main goal, the extended Pruess algorithm would be attractive, but if the main goal is
computation of the eigenvalues, the piecewise constant Pruess method seems a clear winner because of its simplicity
and speed.

\newpage

\begin{table}[h]
\begin{tabular}{cc}
Case 1 & \\
\includegraphics[height= 3in, width= 3in]{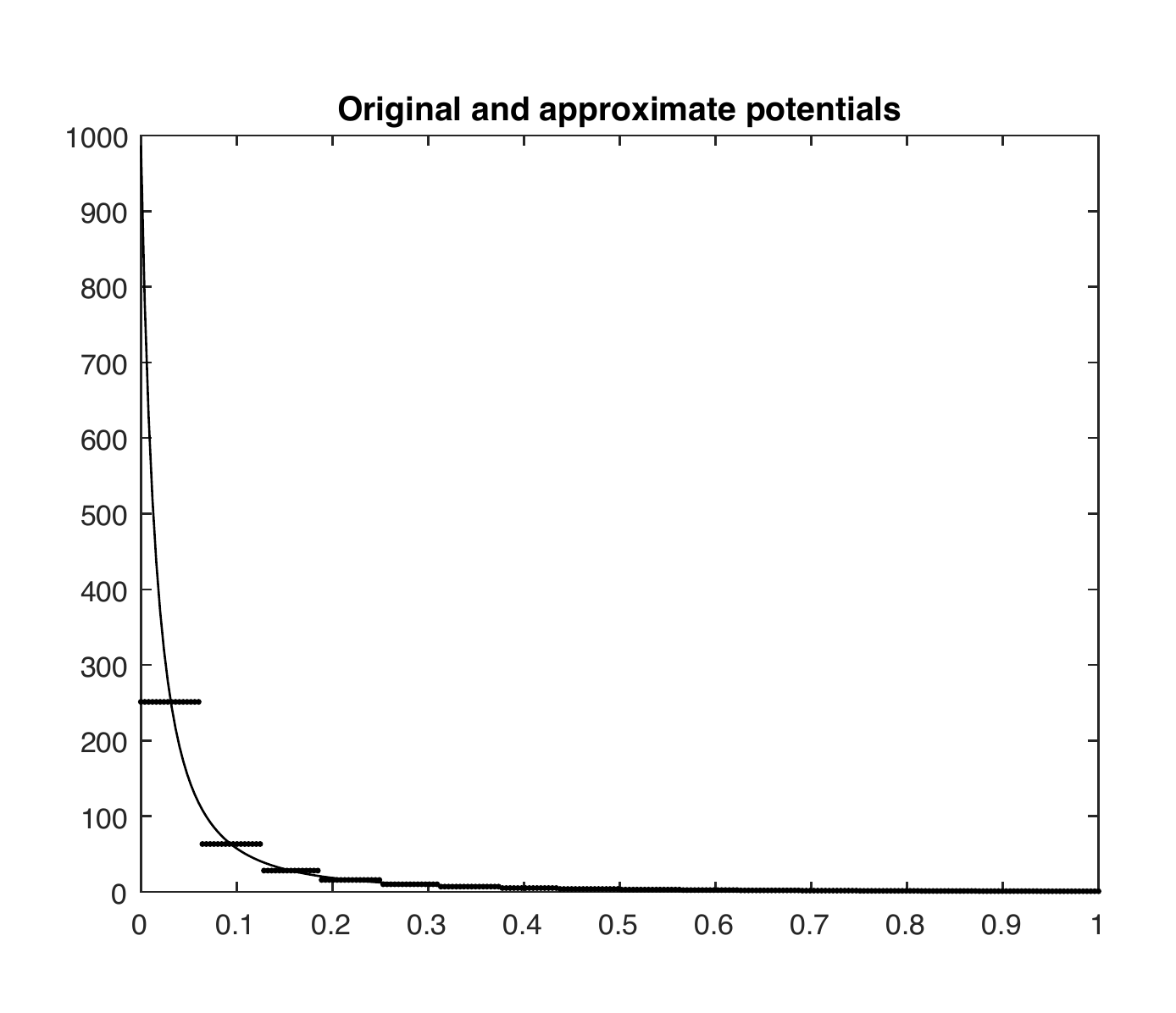}
&
\includegraphics[height= 3in, width= 3in]{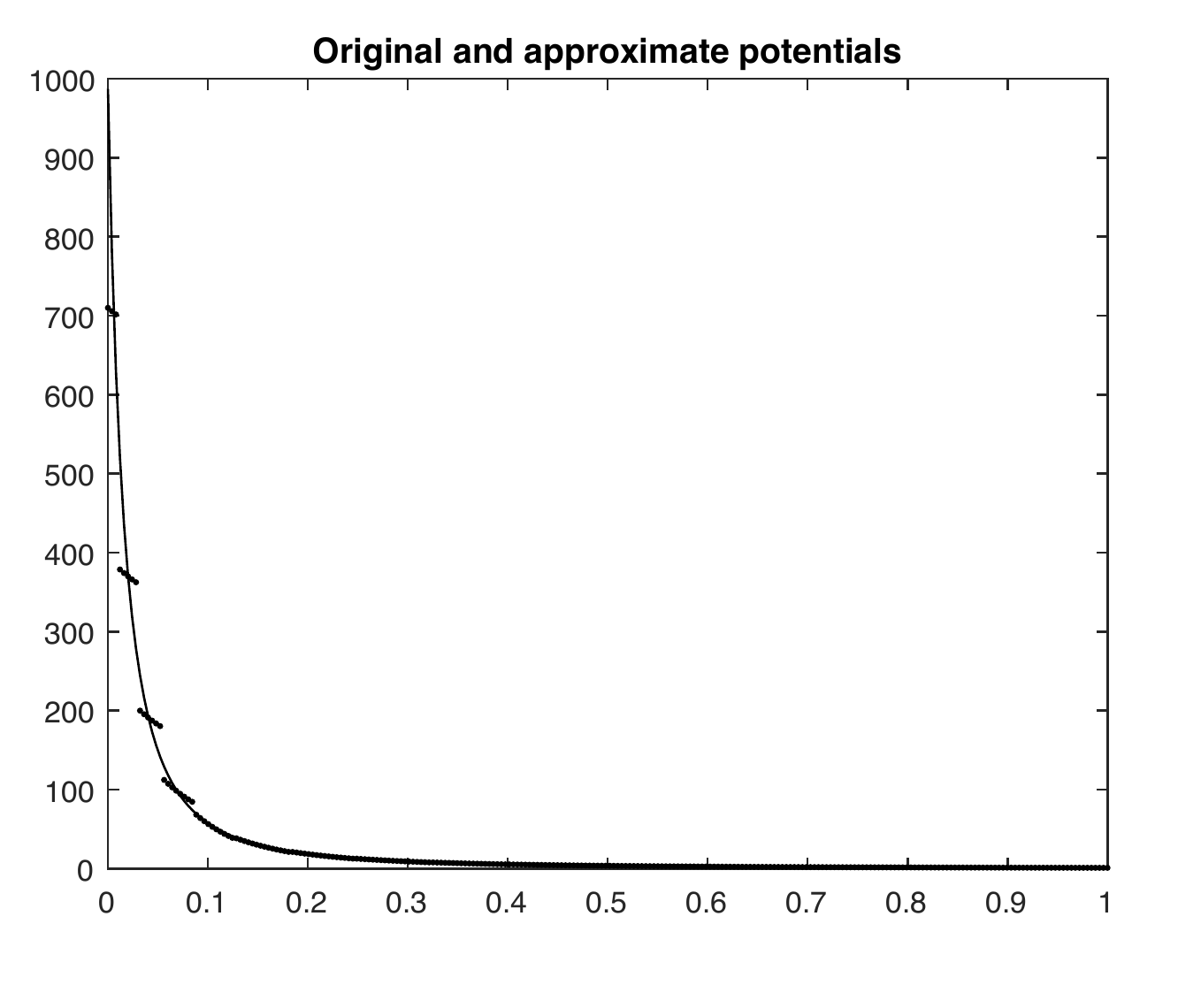} \\
Pruess method, uniform sampling  & 
 Extended method, adaptive sampling
\end{tabular} 
\end{table}

Case 1:
\[ \begin{matrix} & \lambda _1 & \lambda _2 & \lambda _3 & \lambda _{12} & \lambda _{25} \cr
U-P128 & 1.5001 e1 & 4.8792 e1 & 1.0151e2 & 1.4466e3 & 6.1974 e3 \cr
U-P 32 & 1.5015 e1 & 4.8848 e1 & 1.0164e2 & 1.4477e3 & 6.1990e3 \cr
U-P 16 & 1.5055 e1 & 4.9017 e1 & 1.0202e2 & 1.4491e3 & 6.1938e3 \cr
U-X 16 & 1.4938 e1 & 4.8600 e1 & 1.0145e2 & 1.4483e3 & 6.1938e3 \cr
A-P 16 & 1.5031 e1 & 4.8930 e1 & 1.0191e2 & 1.4494e3 & 6.1942 e3 \cr
A-X 16 &  1.4940 e1 & 4.8626 e1 & 1.0127e2 & 1.4463e3 & 6.1969 e3 \cr
best  \ 16 & A-P & A-P & U-X & A-X & A-X
\end{matrix}\]

\newpage

\begin{table}[h]
\begin{tabular}{cc}
Case 2 & \\
\includegraphics[height= 3in, width= 3in]{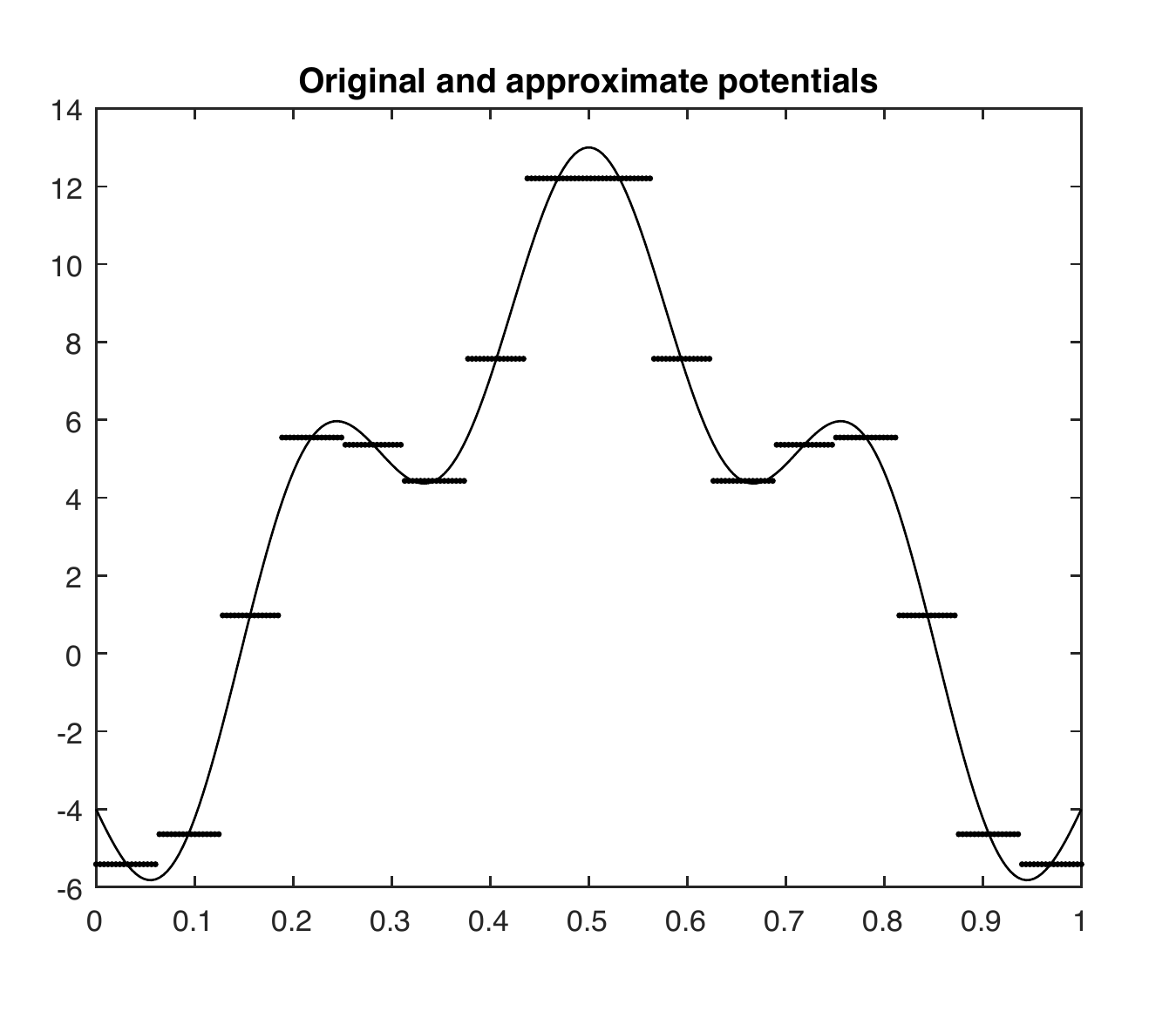}
&
\includegraphics[height= 3in, width= 3in]{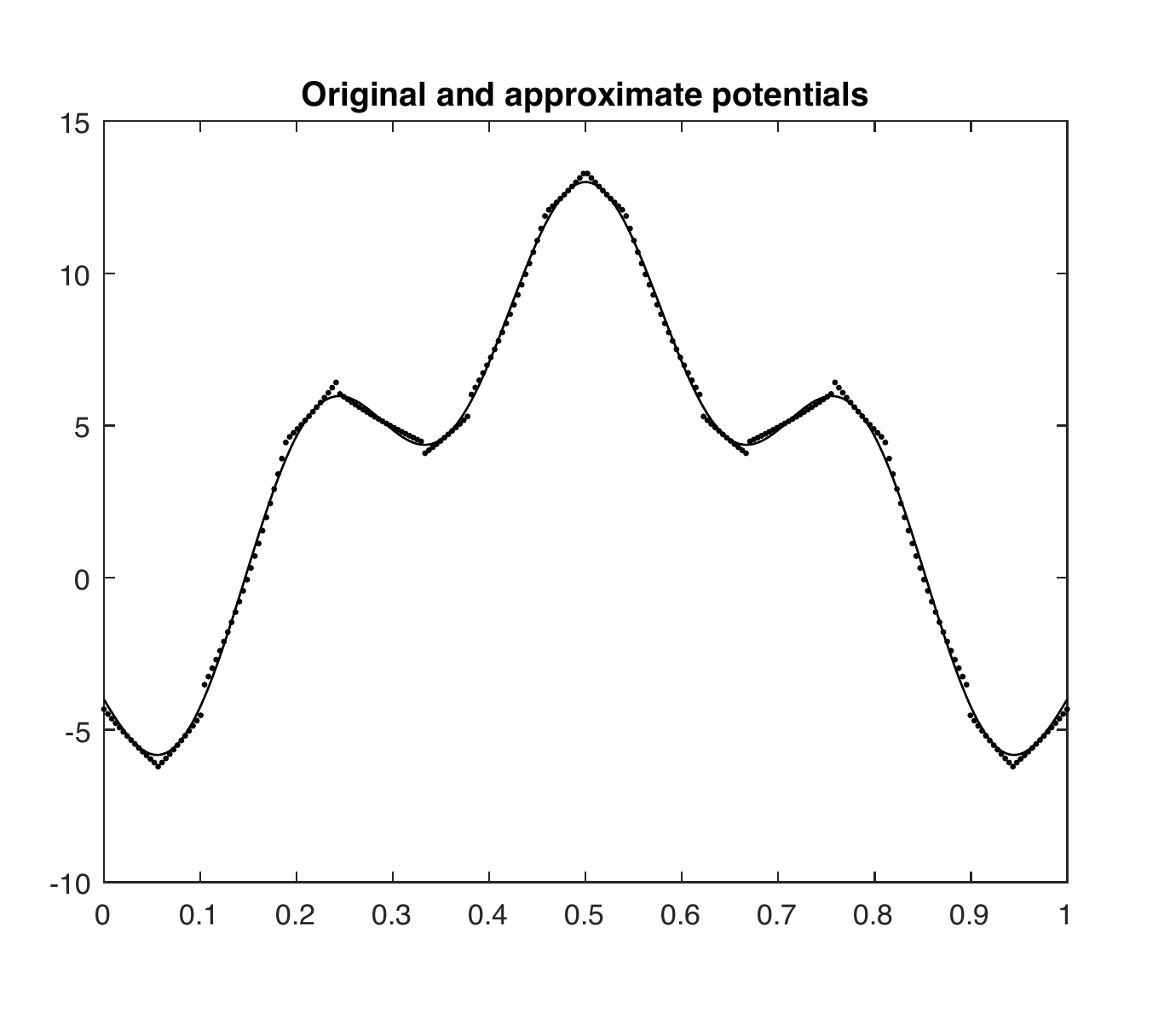} \\
Pruess method, uniform sampling  & 
 Extended method, adaptive sampling
\end{tabular} 
\end{table}

Case 2:
\[ \begin{matrix} & \lambda _1 & \lambda _2 & \lambda _3 & \lambda _{12} & \lambda _{25} \cr
U-P128 & 1.6656e1 & 4.3260e1 & 9.3189e1 & 1.4245 e3 & 6.1718 e3 \cr
U-P 32 & 1.6651e1 & 4.3257e1 & 9.3176 e1 & 1.4245 e3 &  6.1718 e3 \cr
U-P 16 & 1.6635e1 & 4.3247e1 & 9.3135 e1 & 1.4247 e3 &  6.1718 e3 \cr
U-X 16 & 1.6678e1 & 4.3273e1 & 9.3235 e1 & 1.4245 e3 &  6.1718 e3 \cr
A-P 16 & 1.6638e1 & 4.3175e1 & 9.3188 e1 & 1.4245 e3 &  6.1718 e3 \cr
A-X 16 &  1.6664e1 & 4.3271e1 & 9.3199 e1 & 1.4245 e3 &  6.1718 e3 \cr
best \ 16 & A-X & A-X & A-P & - & - 
\end{matrix}\]

\newpage

\begin{table}[h]
\begin{tabular}{cc}
Case 3 & \\
\includegraphics[height= 3in, width= 3in]{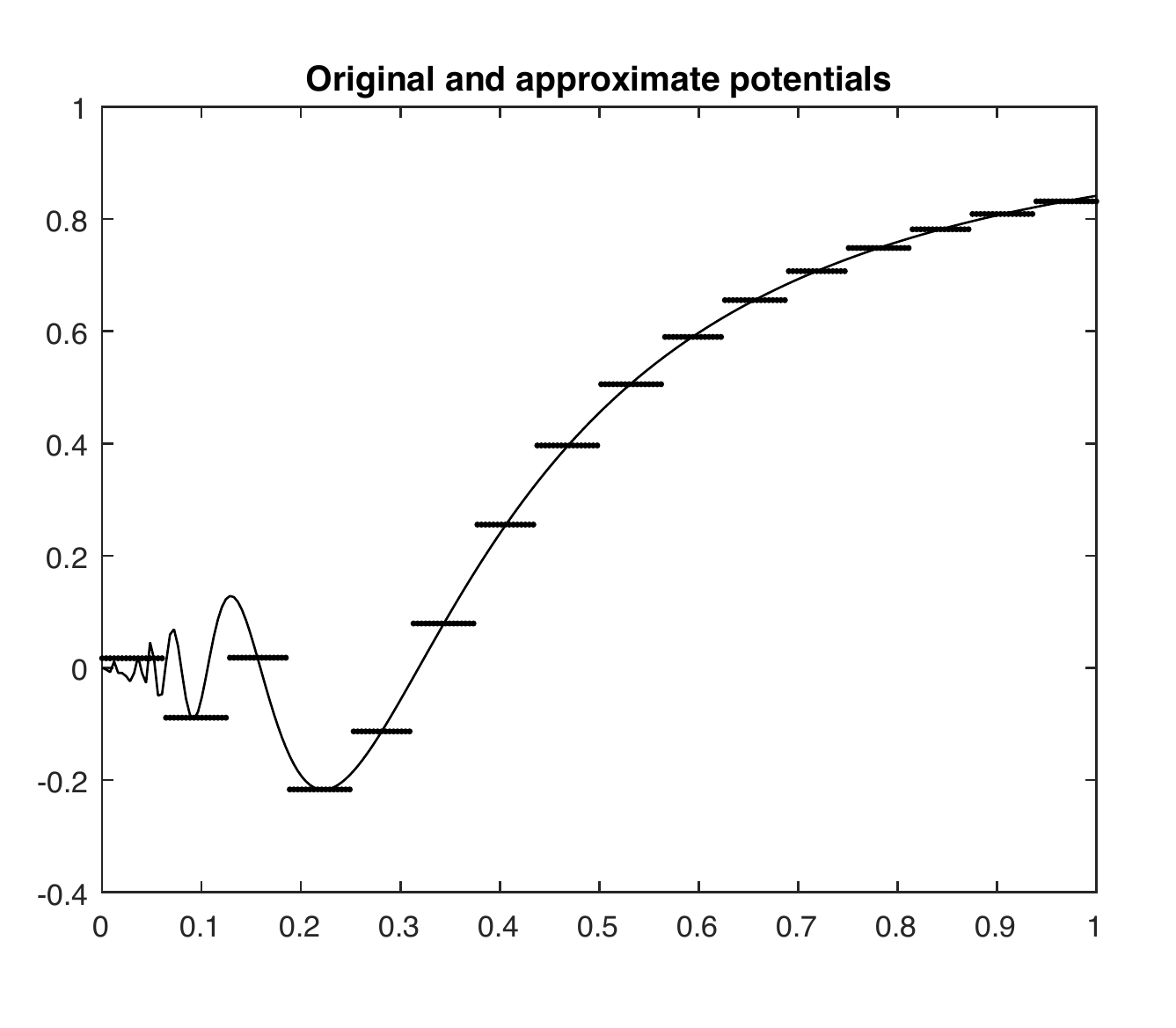}
&
\includegraphics[height= 3in, width= 3in]{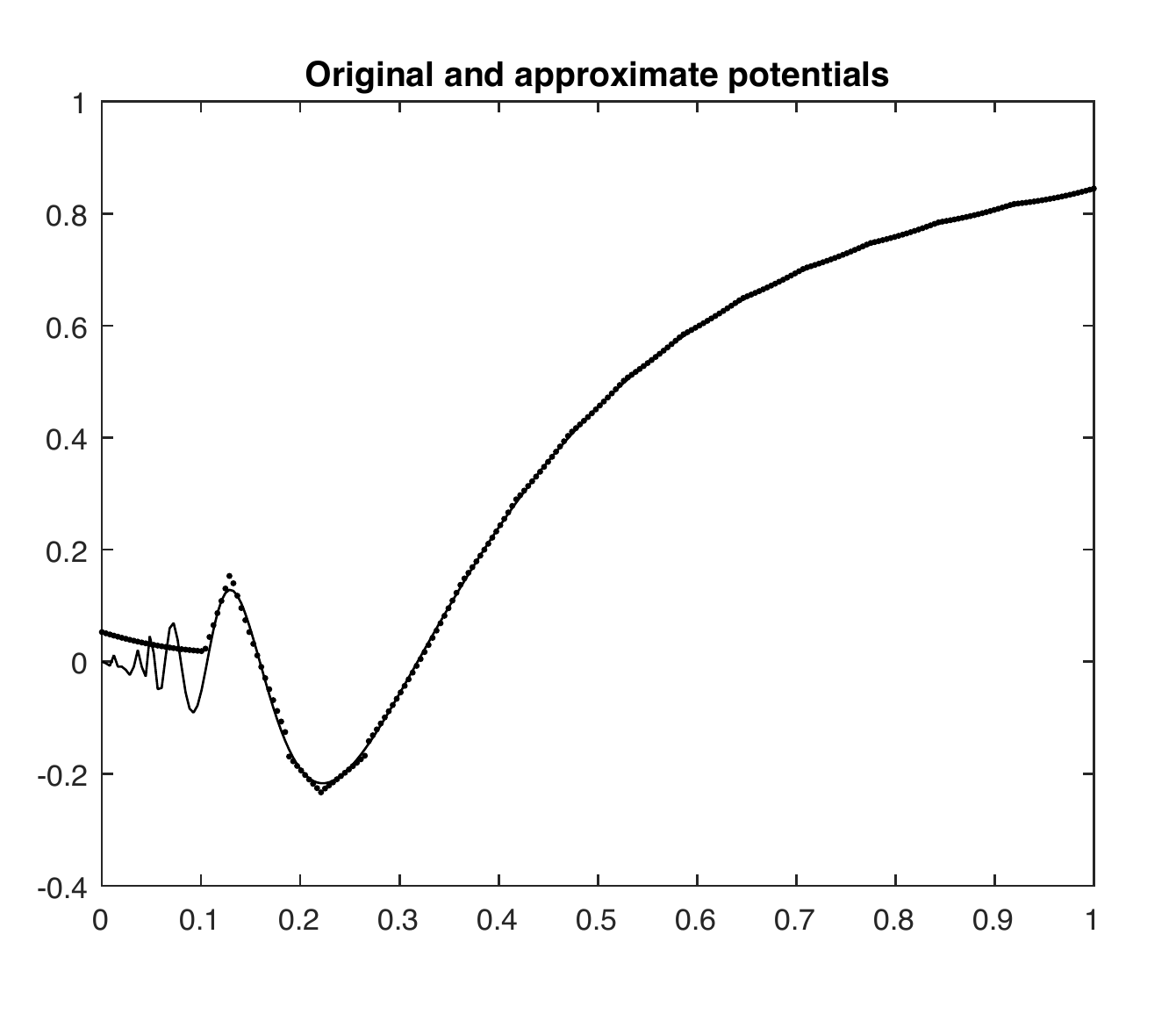} \\
Pruess method, uniform sampling  & 
 Extended method, adaptive sampling
\end{tabular} 
\end{table}

Case 3:
\[\begin{matrix} & \lambda _1 & \lambda _2 & \lambda _3 & \lambda _{12} & \lambda _{25} \cr
U-P128 & 1.0250e1 & 3.9820e1 & 8.9210e1 & 1.4216e3 & 6.1689 e3 \cr
U-P 32 & 1.0250e1 & 3.9821 e1  & 8.9212 e1 & 1.4216 e3 & 6.1689e3 \cr
U-P 16 & 1.0249e1 & 3.9818 e1  & 8.9204 e1 & 1.4216 e3 & 6.1689e3 \cr
U-X 16 & 1.0249e1 & 3.9816 e1  & 8.9204 e1 & 1.4216 e3 & 6.1689e3 \cr
A-P 16 &1.0248e1 & 3.9815e1 & 8.9202e1 & 1.4216e3 & 6.1689 e3 \cr
A-X 16 & 1.0250e1 & 3.9821e1 & 8.9214e1 & 1.4216e3 & 6.1689 e3 \cr
best \ 16 & A-X & A-X & A-X & - & - 
\end{matrix}\]

\newpage

\begin{table}[h]
\begin{tabular}{cc}
Case 4 & \\
\includegraphics[height= 3in, width= 3in]{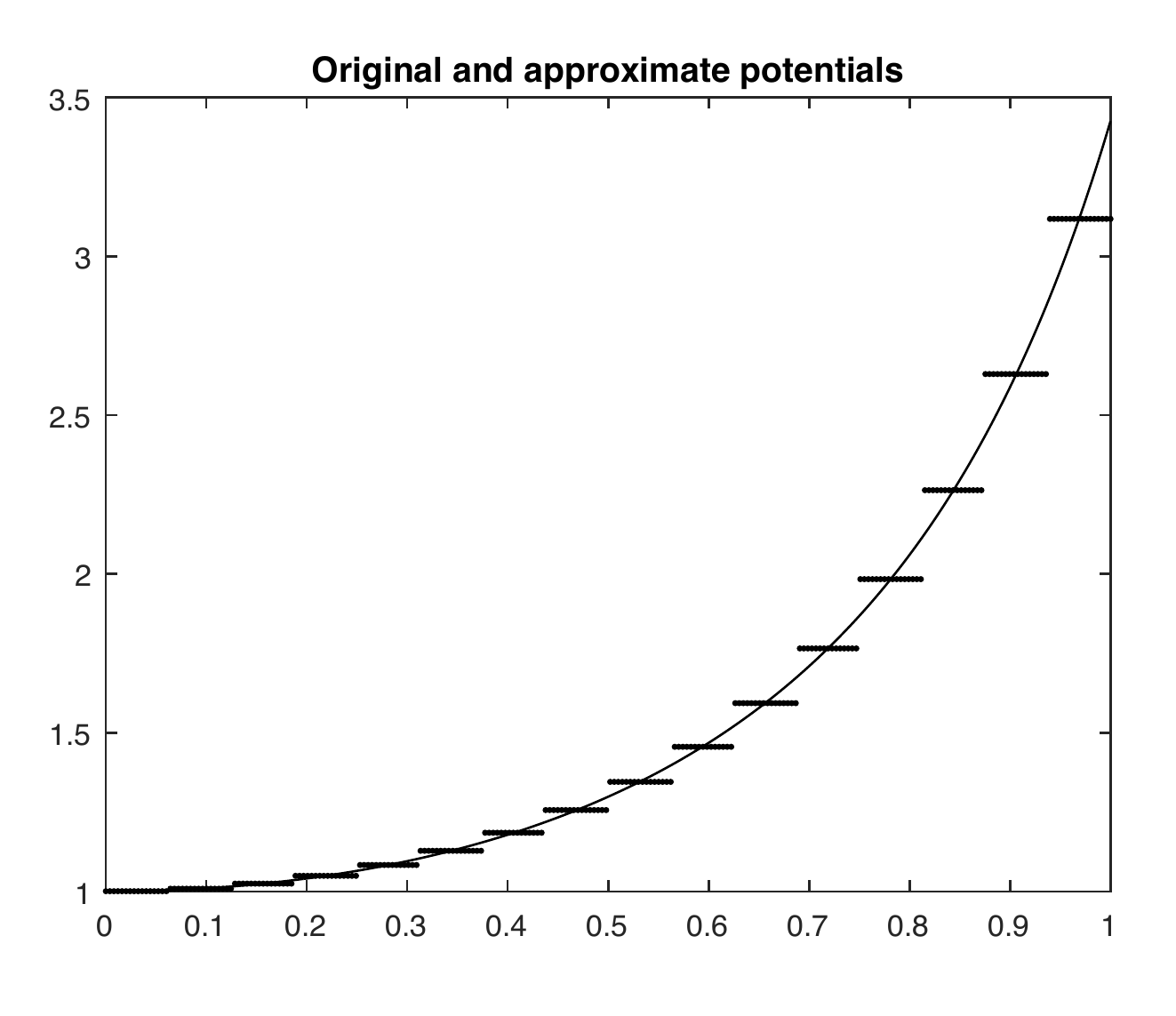}
&
\includegraphics[height= 3in, width= 3in]{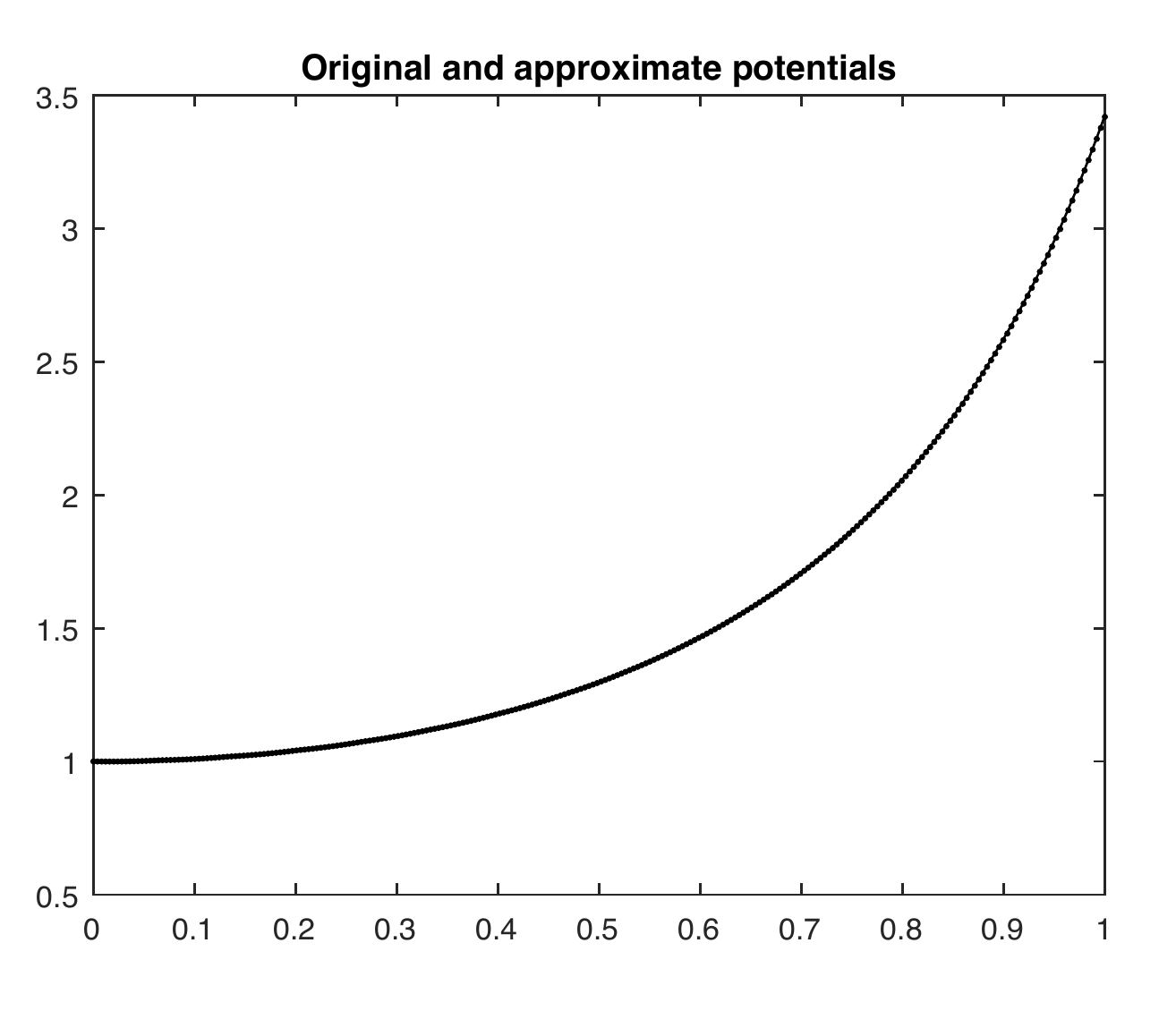} \\
Pruess method, uniform sampling  & 
 Extended method, adaptive sampling
\end{tabular} 
\end{table}

Case 4:
\[ \begin{matrix} & \lambda _1 & \lambda _2 & \lambda _3 & \lambda _{12} & \lambda _{25} \cr
U-P128 & 1.1255e1 & 4.0979e1 & 9.0357e1 & 1.4228e3 & 6.1701e3 \cr
U-P 32 & 1.1256e1& 4.0980e1 & 9.0357e1 & 1.4228e3 & 6.1701e3 \cr
U-P 16 & 1.1256e1& 4.0981e1 & 9.0359e1 & 1.4228e3 & 6.1701e3 \cr
U-X 16 & 1.1254e1& 4.0978e1 & 9.0355e1 & 1.4228e3 & 6.1701e3 \cr
A-P 16 & 1.1256e1 & 4.0980e1 & 9.0357e1 & 1.4228e3 & 6.1701e3 \cr
A-X 16 & 1.1254e1 & 4.0978e1 & 9.0356e1 & 1.4228e3 & 6.1701e3 \cr
best - 16 & - & - & A-P & - & -
\end{matrix} \]

\newpage

\begin{table}[h]
\begin{tabular}{cc}
Case 5 & \\
\includegraphics[height= 3in, width= 3in]{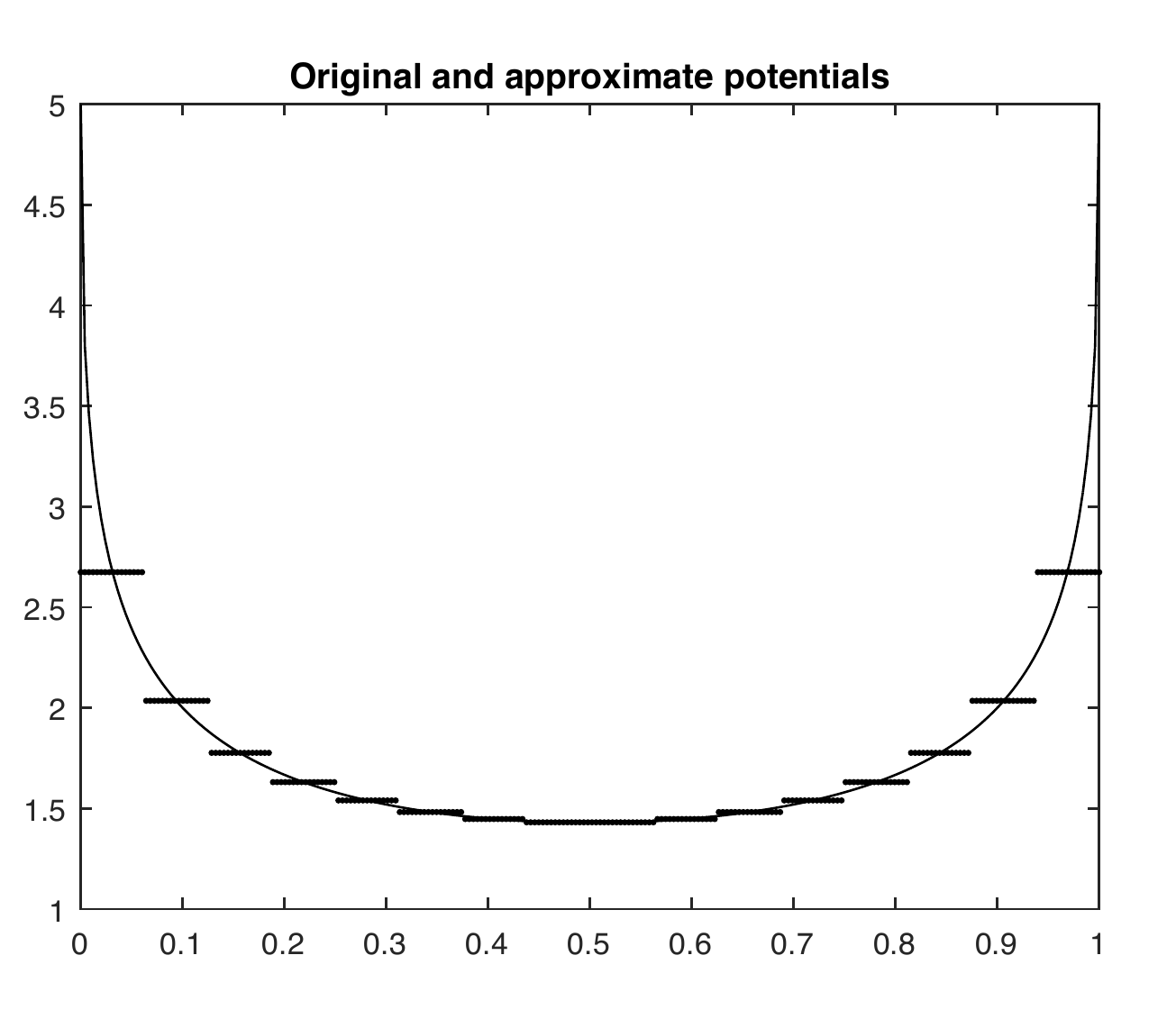}
&
\includegraphics[height= 3in, width= 3in]{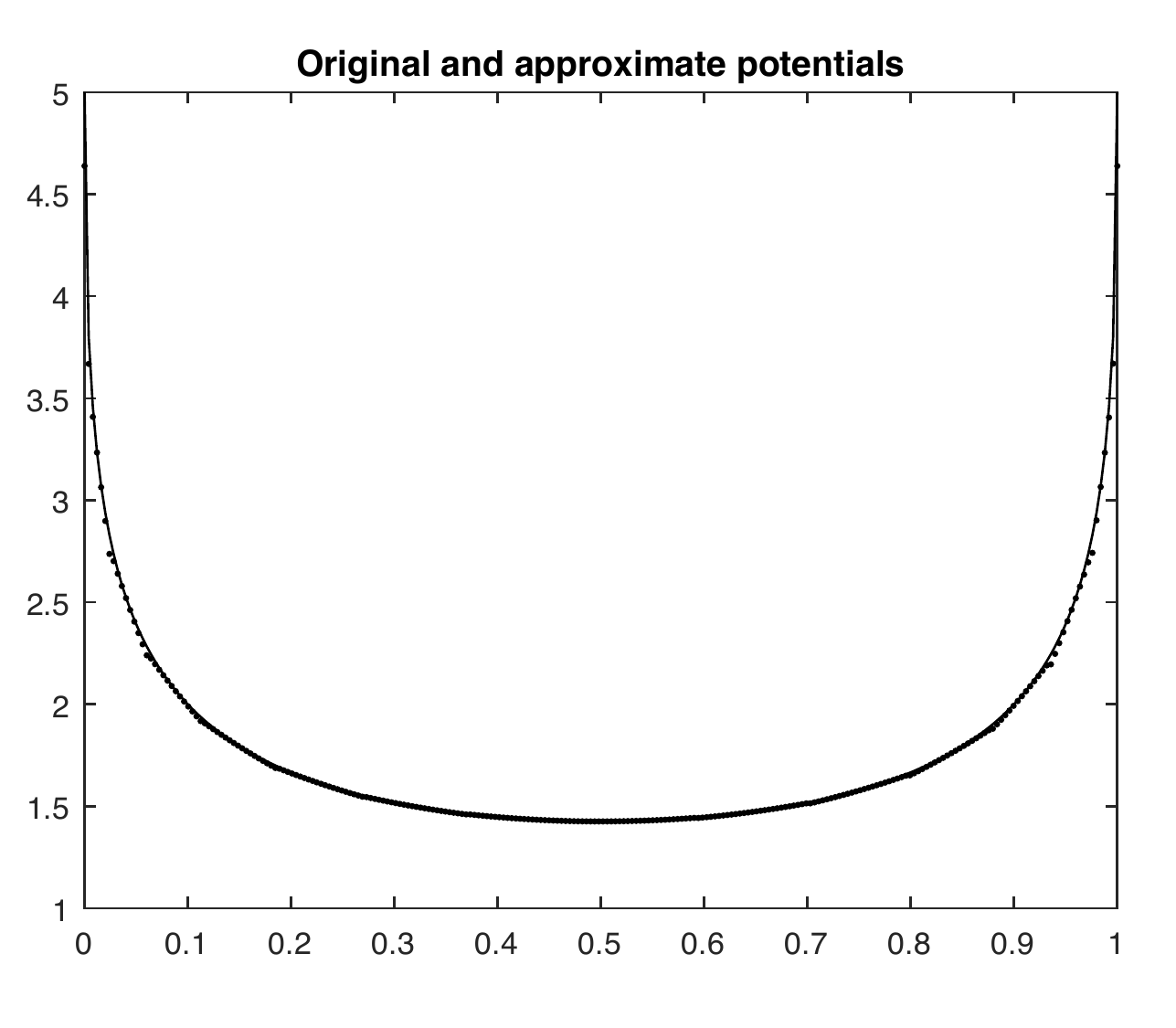} \\
Pruess method, uniform sampling  & 
 Extended method, adaptive sampling
\end{tabular} 
\end{table}

Case 5:
\[ \begin{matrix} & \lambda _1 & \lambda _2 & \lambda _3 & \lambda _{12} & \lambda _{25} \cr
U-P128 & 1.1385e1 & 4.1111e1 & 9.0504e1 & 1.4230e3 & 6.1703e3 \cr
U-P 32 & 1.1385e1 & 4.1111e1 & 9.0506e1 & 1.4230e3 & 6.1703e3 \cr
U-P 16 & 1.1386e1 & 4.1114e1 & 9.0510e1 & 1.4230e3 & 6.1703e3 \cr
U-X 16 & 1.1382e1 & 4.1102e1 & 9.0488e1 & 1.4229e3 & 6.1703e3 \cr
A-P 16 &  1.1384e1 & 4.1108e1 & 9.0504e1 & 1.4230e3 & 6.1703e3 \cr
A-X 16 &  1.1382e1 & 4.1106e1 & 9.0499e1 & 1.4230e3 & 6.1703e3 \cr
best \ 16 & P & P & A-P & - & -
\end{matrix}\]

\newpage

\bibliographystyle{amsalpha}

\end{document}